\documentclass[3p,preprint]{elsarticle_preprint}

\usepackage{amsmath}
\usepackage{amsfonts}
\usepackage{amssymb,latexsym}
\usepackage{amsthm}
\usepackage{graphicx}

\newtheorem{definition}{Definition}[section]
\newtheorem{theorem}{Theorem}[section]

\newtheorem{remark}{Remark}[section]

\usepackage{color}
\usepackage{newlfont}
\usepackage{subfigure}
\usepackage{verbatim}
\usepackage{rotating}
\usepackage{multirow}
\usepackage{units}
\usepackage{bm}
\usepackage[english]{babel}
\usepackage[utf8]{inputenc}
\usepackage{algorithm}
\usepackage[noend]{algpseudocode}
\usepackage{booktabs}
\usepackage{caption}
\usepackage{hyperref}

\newcommand{\RR}{\mathbb{R}}

\newcommand{\EE}{\mathbb{E}}

\newcommand{\mA}{\mathsf{A}}
\newcommand{\mL}{\mathsf{L}}

\journal{R. Cavoretto A. De Rossi}

\begin{document}

\begin{frontmatter}

\title{An adaptive residual sub-sampling algorithm for kernel interpolation based on maximum likelihood estimations}



\author[address-TO,address-GNCS]{Roberto Cavoretto} 
\ead{roberto.cavoretto@unito.it}

\author[address-TO,address-GNCS]{Alessandra De Rossi}
\ead{alessandra.derossi@unito.it}



\address[address-TO]{Department of Mathematics \lq\lq Giuseppe Peano\rq\rq, University of Turin, via Carlo Alberto 10, 10123 Turin, Italy}

\address[address-GNCS]{Member of the INdAM Research group GNCS}

\begin{abstract}

In this paper we propose an enhanced version of the residual sub-sampling method (RSM) in \cite{dri07} for adaptive interpolation by radial basis functions (RBFs). More precisely, we introduce in the context of sub-sampling methods a maximum profile likelihood estimation (MPLE) criterion for the optimal selection of the RBF shape parameter. This choice is completely automatic, provides highly reliable and accurate results for any RBFs, and, unlike the original RSM, guarantees that the RBF interpolant exists uniquely. The efficacy of this new method, called MPLE-RSM, is tested by numerical experiments on some 1D and 2D benchmark target functions.    

\end{abstract}



\end{frontmatter}


\section{Introduction} \label{sec:1}

In \cite{dri07} the residual sub-sampling method (RSM) is proposed. This adaptive scheme is based on radial basis function (RBF) interpolation. It is used to approximate the unknown target function on uniformly distributed points, and then the residual is evaluated at halfway points. The latter are added to the point set when the residual is over a prescribed refinement threshold, whereas they are removed from that set when it is under a predefined coarsening threshold. Further, the interpolating process is characterized by a variable selection of the multiquadric RBF shape parameter. The user assigns an initial value of such parameter, and then it is updated point-by-point based on node spacing. A similar strategy has also been implemented in \cite{zha17}, though in both cases the change of the shape parameter at each point does not guarantee the invertibility of the interpolation matrix. Indeed, as known in literature and observed in \cite{dri07}, the use of node-dependent shape parameters breaks the symmetry of the interpolation matrix as well as the proof of its nonsingularity.  

While several adaptive schemes exist for solving boundary value problems or partial differential equations (see e.g. \cite{cav20b,cav20c,esm12,qia21}), the problem of constructing adaptive algorithms in RBF interpolation has been considered only in part (see \cite{beh02,boz02,dri07,gao20,gao20b,zha17} and references therein). This fact motivates us to enhance further the current state of the art on the topic.

In this paper we present a modification of the original RSM in \cite{dri07}, proposing an optimal selection of the shape parameter via a \emph{maximum profile likelihood estimation} (MPLE) criterion, which relies on a solid stochastic framework (see \cite{fas11,sch11,sch13}). This choice is totally automatic, i.e., user's action is not required, either initially, but a single (near-optimal) parameter is determined for every node. In this way, the interpolation problem is well-posed and hence the RBF (or kernel) interpolant exists uniquely, naturally provided that the kernel matrix is positive definite (see e.g. \cite{fas15}). Moreover, the use of MPLE technique is particularly useful to take under control the ill-conditioning of the interpolation matrix when in the iterative/adaptive method the number of interpolation points grows. Therefore, unlike the previous methods discussed in \cite{dri07,zha17}, this modification of the RSM, called MPLE-RSM, provides a twofold benefit. As a matter of fact, it enables to solve the above mentioned issue of nonsingularity of the interpolation matrix, and give \lq\lq good\rq\rq\ predictions of the shape parameter for any kernel avoiding user's tuning each time. As our numerical results show, application of the MPLE criterion generally results in an adaptive interpolation scheme more efficient than the original RSM, also reducing the number of points needed for satisfying the expected precision. The improved method is tested by taking some 1D and 2D benchmark target functions. 

The paper is organized as follows. In Section \ref{sec:2} we introduce the kernel based interpolation. Section \ref{sec:3} presents the MPLE strategy to determine the optimal shape parameter in the improved RSM. In Section \ref{sec:4} we describe the adaptive algorithm. In Section \ref{sec:5} we show some numerical results, illustrating the performance of the MPLE-RSM and comparing the latter with the basic RSM discussed in \cite{dri07}. Section \ref{sec:6} contains conclusions.

\section{Kernel based interpolation} \label{sec:2}
Kernel based methods are powerful tools for data interpolation. In this section we introduce basic notations and results for kernel based interpolation. For further theoretical background and other details, we refer the reader to \cite{buh03,fas15,wen05}.	

Given a compact domain $\Omega \subset \RR^d$, we assume that the $N$ distinct data points (or nodes) are defined by the set $X = \{\boldsymbol{x}_i\}_{i=1}^{N} \subseteq \Omega$. The data values associated with $\boldsymbol{x}_i$, $i=1,\ldots,N$, are given by $y_i = f(\boldsymbol{x}_i) \in \RR$, where the latter are obtained by sampling some function $f: \Omega \rightarrow \mathbb{R}$. Thus, we want to find a function $s_X : \Omega \rightarrow \mathbb{R}$ satisfying the interpolation conditions
\begin{align} \label{ic}
s_X(\boldsymbol{x}_i) = y_i, \qquad i =1,\ldots,N.
\end{align}

We express the interpolant $s_X$ in terms of a kernel $\Phi: \Omega \times \Omega \rightarrow \mathbb{R}$, i.e.
\begin{equation} \label{rad2}
	s_X( \boldsymbol{x})= \sum_{j=1}^{N} c_j \Phi (\boldsymbol{x},\boldsymbol{x}_j), \qquad \boldsymbol{x} \in \Omega.
\end{equation} 

If the kernel $\Phi$ is symmetric and strictly positive definite (SPD), the interpolation matrix $\mA = (\mA_{ij})$ with the entries $\mA_{ij}= \Phi(\boldsymbol{x}_i,\boldsymbol{x}_j)$, $i,j=1,\ldots,N$, is positive definite for any set $X$. The coefficients $c_j$ in \eqref{rad2} are uniquely determined by enforcing the interpolation conditions \eqref{ic} and can be obtained by solving the symmetric linear system
\begin{align} \label{linsys}
\mA \boldsymbol{c} = \boldsymbol{y},
\end{align}
where $\boldsymbol{c} = (c_1,\ldots,c_N)^T$ and $\boldsymbol{y} = (y_1,\ldots,y_N)^T$.

Associated with the kernel $\Phi$ in \eqref{rad2} we may define a SPD RBF $\phi : \RR_0^+ \rightarrow \RR$ such that
\begin{align*}
\Phi(\boldsymbol{x},\boldsymbol{x}_j)=\phi_{\varepsilon}(||\boldsymbol{x}-\boldsymbol{x}_j||_2)=\phi(\varepsilon||\boldsymbol{x}-\boldsymbol{x}_j||_2), \qquad \forall \boldsymbol{x},\boldsymbol{x}_j\in \Omega,
\end{align*}
where $\varepsilon > 0$ is the so-called \emph{shape parameter}, and $||\cdot||_2$ denotes the Euclidean norm on $\RR^d$. Moreover, we know that the choice of a \lq\lq good\rq\rq\ value of $\varepsilon$ is generally a crucial task for kernel based interpolation, but at the same time also a big issue (see e.g. \cite{cav21c,gol15}, or \cite[Chapter 14]{fas15}). Some examples of popular SPD RBFs (or radial kernels) together with their smoothness degrees and abbreviations are listed as follows (see \cite{fas07, wen05}):
\begin{align*}
\phi_{\varepsilon}(r) =\left\{
\begin{array}{llllc}
\exp(-\varepsilon^2 r^2), & & \quad \mbox{Gaussian $C^{\infty}$}, & & \quad \mbox{GA}  \medskip  \\ 
(1+\varepsilon^2r^2)^{-1/2}, & & \quad \mbox{Inverse MultiQuadric $C^{\infty}$}, & & \quad  \mbox{IMQ}  \medskip  \\
\exp(-\varepsilon r) (\varepsilon^3 r^3 + 6\varepsilon^2r^2+15\varepsilon r+15),  & &	\quad \mbox{Mat$\acute{\text{e}}$rn $C^6$}, & & \quad  \mbox{M6}  \medskip  \\
\exp(-\varepsilon r) (\varepsilon^2r^2+3\varepsilon r+3), & & \quad \mbox{Mat$\acute{\text{e}}$rn $C^4$}, & & \quad \mbox{M4}  \medskip  \\
\exp(-\varepsilon r) (\varepsilon r+1), & & \quad \mbox{Mat$\acute{\text{e}}$rn $C^2$}, & & \quad \mbox{M2}  
\end{array}
\right.
\end{align*}

When solving the linear system \eqref{linsys}, the solution is often very sensitive to changes in the data. Moreover, such sensitivity is influenced by the choice of the shape parameter $\varepsilon$. A criterion for measuring the numerical stability of a kernel method is to compute the condition number of the interpolation matrix $\mA$. Hence, since the kernel $\Phi$ is symmetric and SPD, the condition number of $\mA$ is defined as follows:
\begin{align} \label{condnum}
\kappa(\mA) = ||\mA||_2 ||\mA^{-1}||_2 = \frac{\lambda_{\max}}{\lambda_{\min}},
\end{align}
where $\lambda_{\max}$ and $\lambda_{\min}$ are the largest and smallest eigenvalues of $\mA$.

Furthermore, for the kernel $\Phi$ there exists the so-called \emph{native space}, which is a Hilbert space ${\cal N}_{\Phi}(\Omega)$ 
with inner product $(\cdot,\cdot)_{{\cal N}_{\Phi}(\Omega)}$ in which the kernel $\Phi$ is reproducing, i.e., for any $f \in {\cal N}_{\Phi}(\Omega)$ we have the identity $f(\boldsymbol{x}) = (f,\Phi(\cdot,\boldsymbol{x}))_{{\cal N}_{\Phi}(\Omega)}$, with $\boldsymbol{x}\in \Omega$. Then, if we introduce a pre-Hilbert space $H_{\Phi}(\Omega)= \mbox{span}\{\Phi(\cdot,\boldsymbol{x}),$ $\boldsymbol{x} \in \Omega\}$, with reproducing kernel $\Phi$ and equipped with the bilinear form $(\cdot,\cdot)_{H_{\Phi}(\Omega)}$, the native space ${\cal N}_{\Phi}(\Omega)$ of $\Phi$ is its completion with respect to the norm $||\cdot||_{H_{\Phi}(\Omega)}=\sqrt{(\cdot,\cdot)_{H_{\Phi}(\Omega)}}$. In particular, for all $f \in {H_{\Phi}(\Omega)}$ we have $||f||_{{\cal N}_{\Phi}(\Omega)}=||f||_{H_{\Phi}(\Omega)}$ (see \cite{wen05}). Now, we can thus provide an error bound in terms of the well-known \textsl{power function} $P_{\Phi, X}$  (see e.g. \cite[Theorem 14.2]{fas07}):
\begin{theorem} \label{th1}
Let $\Omega \subseteq \RR^d$, $\Phi \in C(\Omega \times \Omega)$ be strictly positive definite on $\RR^d$, and suppose that $X=\{\boldsymbol{x}_i\}_{i=1}^N$ has distinct points. Then, for all $f\in {\cal N}_{\Phi}(\Omega)$, we have
\begin{align*} 
	|f(\boldsymbol{x})-s_X(\boldsymbol{x})| \leq P_{\Phi,X}(\boldsymbol{x}) ||f||_{{\cal N}_{\Phi(\Omega)}}, \quad \boldsymbol{x} \in \Omega.
\end{align*}
\end{theorem}
The generic error estimate of Theorem \ref{th1} can further be refined as shown in \cite[Theorem 14.5]{fas07}:
\begin{theorem} \label{th2}
	Let $\Omega\subseteq \RR^d$ be bounded and satisfy an interior cone condition. Suppose that $\Phi \in C^{2k}(\Omega\times \Omega)$ is symmetric and strictly positive definite. Then, for all $f\in {\cal N}_{\Phi}(\Omega)$, there exist constants $h_0$, $C > 0$ (independent of $\boldsymbol{x}$, $f$ and $\Phi$) such that
\begin{align*}
	\left|f(\boldsymbol{x})-s_X(\boldsymbol{x})\right|\leq C h_{X,\Omega}^{k} \sqrt{C_{\Phi}(\boldsymbol{x})} \left\|f\right\|_{{\cal N}_{\Phi}(\Omega)},
\end{align*}
provided $h_{X,\Omega}\leq h_0$. Here
\begin{align*}
   C_{\Phi}(\boldsymbol{x})= \max_{\left|\boldsymbol{\beta}\right|=2k,} \max_{\boldsymbol{w},\boldsymbol{z}\in \Omega\cap B(\boldsymbol{x},c_2 h_{X,\Omega})} \left|D_{2}^{\boldsymbol{\beta}}\Phi(\boldsymbol{w},\boldsymbol{z})\right|
\end{align*}
with $B(\boldsymbol{x},c_2 h_{X,\Omega})$ denoting the ball of radius $c_2 h_{X,\Omega}$ centred at $\boldsymbol{x}$, and $h_{X, \Omega}$ being the \emph{fill distance} 
\begin{align*} 
	h_{X, \Omega} =  \sup_{ \boldsymbol{x} \in \Omega} \min_{ \boldsymbol{x}_j  \in X} || \boldsymbol{x} - \boldsymbol{x}_j||_2.
\end{align*}
\end{theorem}

Theorem \ref{th2} states that interpolation with a $C^{2k}$ smooth kernel $\Phi$ has approximation order $k$. Thus, we deduce that: (i) for $C^{\infty}$ SPD kernels, the approximation order $k$ is arbitrarily high; (ii) for SPD kernels with limited smoothness, the approximation order is limited by the smoothness of the kernel. For more refined error estimates, we refer the reader to the monograph \cite{wen05}.

\section{MPLE criterion for near-optimal choice of the shape parameter} \label{sec:3}

In Section \ref{sec:2} we compute the interpolant $s_X$ in \eqref{rad2} by solving the linear system \eqref{linsys}, where the kernel matrix $\mA$ is symmetric and positive definite. However, by the \emph{uncertainty} or \emph{trade-off} \emph{principle} \cite{sch95} we know that using a standard RBF one cannot have high accuracy and stability at the same time. In fact, when the best level of accuracy is typically achieved, i.e., in the \emph{flat limit} $\varepsilon \rightarrow 0$, the interpolation matrix may be very ill-conditioned. It is therefore important to study a criterion that enables us to make reliable $\varepsilon$-predictions. In this work we discuss the MPLE, which we will apply in the residual sub-sampling interpolation method. 

\subsection{Gaussian random field and density function} \label{sec:31}
	
The MPLE criterion is mainly based on a stochastic framework, and so the concept of Gaussian random field (or Gaussian process) is introduced \cite{fas15}. 

\begin{definition}
The random field $Y=\{Y_{\boldsymbol{x}}\in \Omega\}$ is called a Gaussian random field if, for any given choice of finitely many distinct points $X=\{\boldsymbol{x}_i\}_{i=1}^{N} \subseteq \Omega$, the vector of random variable $\boldsymbol{Y}=(Y_{\boldsymbol{x}_1},\ldots,Y_{\boldsymbol{x}_N})^T$ has a multivariate normal distribution with mean vector $\boldsymbol{\mu} = \EE[\boldsymbol{Y}]$ and covariance matrix $\sigma^2\mA=(Cov(Y_{\boldsymbol{x}_i},Y_{\boldsymbol{x}_j}))_{i,j=1}^N$, where $\sigma^2$ is the process variance.
\end{definition}

In terms of notation we write $\boldsymbol{Y}\sim {\cal N}(\boldsymbol{\mu},\sigma^2\mA)$ to denote that $\boldsymbol{Y}$ is a vector of Gaussian random variables, or $Y\sim {\cal N}(\mu,\sigma^2\mA)$ to indicate that $Y$ is a Gaussian random field.

The \emph{multivariate normal distribution} has the density function
\begin{align} \label{MNDDF}
p_{\boldsymbol{Y}}(\boldsymbol{y}) = \frac{1}{\sqrt{(2\pi\sigma^2)^N\det \mA}}\exp\left[-\frac{1}{2\sigma^2}(\boldsymbol{y}-\boldsymbol{\mu})^T\mA^{-1}(\boldsymbol{y}-\boldsymbol{\mu})\right].
\end{align}

In the stochastic setting the process variance plays an important role, for instance, in the formulation of the kriging variance and in parameter estimation. In fact, it does not affect the kernel interpolant (or kriging predictor, as known in this context), but this influences its variance and as a consequence the maximum likelihood estimation in the choice of the optimal value of $\varepsilon$. 

Now, the parameters $\varepsilon$ and $\sigma^2$ might be viewed as draws from random variables ${\cal E}$ and $\Sigma$, respectively, with unknown distributions. By studying the joint distribution $(\Sigma,{\cal E},\boldsymbol{Y})$, the kernel parametrization would require to optimize for both $\sigma^2$ and $\varepsilon$ maximizing $p_{\Sigma,{\cal E} | \boldsymbol{Y}}(\sigma^2,\varepsilon|\boldsymbol{Y}=\boldsymbol{y})$, which we can suppose proportional to the density function 
\begin{align} \label{DF}
p_{\boldsymbol{Y}|\Sigma,{\cal E}}(\boldsymbol{y}|\Sigma = \sigma^2,{\cal E}=\varepsilon) = 
\frac{1}{\sqrt{(2\pi\sigma^2)^N\det \mA}}\exp\left[-\frac{1}{2\sigma^2}\boldsymbol{y}^T\mA^{-1}\boldsymbol{y}\right].
\end{align}
where though $\varepsilon$ does not explicitly appear on the right hand side, it appears within $\mA$.

Note that the function \eqref{DF} derives from \eqref{MNDDF} by using an appropriate notation and setting $\boldsymbol{\mu}=\boldsymbol{0}$. Thus, if we assume that $\boldsymbol{\mu} = \boldsymbol{0}$, then the kernel based interpolant is defined by the linear system \eqref{linsys}. In addition, since $\Phi$ is a SPD kernel, the matrix $\mA$ is positive definite, and so invertible.

\subsection{Determination of the MPLE criterion for kernel interpolation} \label{sec:32}

While the discussion given in Subsection \ref{sec:31} would result in a two-dimensional optimization problem, here we use another technique known as \emph{profile likelihood} in which we define $\sigma^2$ as a function of $\varepsilon$, i.e., $\sigma^2=\sigma^2(\varepsilon)$. Thus, our goal reduces to finding an optimal process variance $\sigma_{\mbox{opt}}^2$ by maximizing $p_{\Sigma | {\cal E}, \boldsymbol{Y}}(\sigma^2 | {\cal E} = \varepsilon, \boldsymbol{Y}=\boldsymbol{y}) \propto p_{\boldsymbol{Y}|\Sigma, {\cal E}}(\boldsymbol{y} | \Sigma =\sigma^2, {\cal E} = \varepsilon)$, see \cite{fas15}.

The concept of maximizing the likelihood function requires the maximization of $p_{\boldsymbol{Y}|\Sigma, {\cal E}}(\boldsymbol{y} | \Sigma =\sigma^2, {\cal E} = \varepsilon)$. However, the optimal value of $\sigma^2$ can be determined by minimizing the negative logarithm of \eqref{DF} (multiplying by 2), i.e., 
\begin{align} \label{neglog}
- 2\log\left(p_{\boldsymbol{Y}|\Sigma, {\cal E}}(\boldsymbol{y} | \Sigma =\sigma^2, {\cal E} = \varepsilon)\right) = N\log 2\pi + N\log \sigma^2 +  \log \det \mA + \frac{1}{\sigma^2} \boldsymbol{y}^T \mA^{-1} \boldsymbol{y}.
\end{align}
Differentiating \eqref{neglog} w.r.t. $\sigma^2$ and equating to zero, we obtain the optimal profile variance
\begin{align} \label{OPV}
\sigma^2_{\mbox{opt}} = \frac{1}{N} \boldsymbol{y}^T \mA^{-1} \boldsymbol{y}.
\end{align}
Hence, by setting \eqref{OPV} in \eqref{neglog}, the minimization process that involves the profile likelihood gives
\begin{align*} 
- 2\log\left(p_{\boldsymbol{Y}|\Sigma, {\cal E}}(\boldsymbol{y} | \Sigma =\sigma^2_{\mbox{opt}}, {\cal E} = \varepsilon)\right) = N\log\left(\boldsymbol{y}^T \mA^{-1} \boldsymbol{y}\right) +  \log \det \mA + N (1+ \log 2\pi -\log N).
\end{align*}
Now, ignoring the constant term $N (1+ \log 2\pi -\log N)$, the cost function to minimize via the MPLE criterion is
\begin{align} \label{mple_func} 
\mbox{MPLE}(\varepsilon) = N\log\left(\boldsymbol{y}^T \mA^{-1} \boldsymbol{y}\right) +  \log \det \mA.
\end{align}

\begin{remark} \label{rem}
The computation of the MPLE criterion \eqref{mple_func} for a range of $\varepsilon$ values is carried out by applying the Cholesky factorization to the matrix, i.e., $\mA=\mL\mL^T$. In practice, the use of such a factorization simplifies the determinant computation, because in this case $\log(\det(\mA))=\log(\det(\mL\mL^T))=2\log(\det(\mL))=2\sum_{i=1}^N\log(\sigma_i^{\mL})$, where $\sigma_i^{\mL}$ denotes the eigenvalues of $\mL$. Finally, in order to quickly find the optimal value of $\varepsilon$, the minimum of the cost function \eqref{mple_func} can be determined by the \textsc{Matlab} \texttt{fminbnd} function (or, in case, any other minimization routine). 
\end{remark}

\section{Adaptive algorithm based on refinement and coarsening processes} \label{sec:4}

In this section we describe our adaptive algorithm, which is based on a computational procedure. The latter enables us to refine and coarsen the distribution of interpolation points.

\subsection{Residual sub-sampling procedure} \label{sec:41}

First of all, we introduce a sequence of point sets $X^{(0)}$, $X^{(1)}$, $\cdots$, such that $X^{(k+1)}$ is generated from $X^{(k)}=\{\boldsymbol{x}_i^{(k)}\}_{i=1}^{N^{(k)}}$ after applying some refinement and/or coarsening strategies. These updates depend on residual evaluations, which lead to an adaptive residual sub-sampling method. Therefore, the resulting process follows the common paradigm to solve, estimate and refine/coarsen till a criterion stop is satisfied. 

Now, defining a \emph{check} or \emph{test} set $T^{(k)}=\{\boldsymbol{t}_i^{(k)}\}_{i=1}^{{N_{T^{(k)}}}} \subset \Omega$, for $k\geq 0$, we can evaluate the residual
\begin{align} \label{resid}
\xi(\boldsymbol{t}_i^{(k)}) = \left|s_{X^{(k)}}(\boldsymbol{t}_i^{(k)}) - f(\boldsymbol{t}_i^{(k)})\right|, \qquad \boldsymbol{t}_i^{(k)} \in T^{(k)},
\end{align}
where $s_{X^{(k)}}$ is the interpolating function defined on the set ${X^{(k)}}$, $N^{(k)}$ being the number of points in $T^{(k)}$. 

The residual error defined in \eqref{resid} measures the deviation between the approximate solution and the function value at the point $\boldsymbol{t}_i^{(k)}$. Thus, when $\boldsymbol{t}_i^{(k)}$ lies in a smooth region, the absolute error $\xi(\boldsymbol{t}_i^{(k)})$ is expected to be small, whereas in the region of less regularity for $f$, or around discontinuities, the residual error $\xi(\boldsymbol{t}_i^{(k)})$ is expected to be large. Notice that for $k=0$ the check set $T^{(0)}$ is defined by starting from $X^{(0)}$, while for $k\geq 1$ the check set $T^{(k)}$ is dependent from $X^{(k)}$ and $X^{(k-1)}$.

Thus, the residual $\xi(\boldsymbol{t}_i^{(k)})$ in \eqref{resid} is used as a criterion to define a refinement set $X^{(k)}_{\mbox{refine}}$ and a coarsening set $X^{(k)}_{\mbox{coarse}}$. In doing so, we need to introduce two tolerances (or thresholds) $\theta_{\mbox{refine}}$ and $\theta_{\mbox{coarse}}$, such that $0 < \theta_{\mbox{coarse}} < \theta_{\mbox{refine}}$. When the value of \eqref{resid} is larger than $\theta_{\mbox{refine}}$, we add the point $\boldsymbol{t}_i^{(k)}$ in the refinement set $X^{(k)}_{\mbox{refine}}$, and so at next step the set $X^{(k)}$ needs to be replaced by $X^{(k)} \cup X^{(k)}_{\mbox{refine}}$. Instead, whenever the error $\xi(\boldsymbol{t}_i^{(k)})$ is smaller than $\theta_{\mbox{coarse}}$, we move a point from the active node set $X^{(k)}$ into the coarsening set $X^{(k)}_{\mbox{coarse}}$, and so $X^{(k)}$ is then updated with $X^{(k)} \backslash X^{(k)}_{\mbox{coarse}}$. As a consequence, at ($k+1$)-step of our adaptive process the set $X^{(k)}$ is updated by adding the refinement set $X^{(k)}_{\mbox{refine}}$ and deleting the coarsening set $X^{(k)}_{\mbox{coarse}}$, that is $X^{(k+1)} = \left\{ X^{(k)} \cup X^{(k)}_{\mbox{refine}} \right\} \backslash X^{(k)}_{\mbox{coarse}}$. The iterative method concludes once the process of addition and/or removal was completed, returning the final set $X^{(k^*)}$, where $k^*$ denotes the last iteration. A pseudo-code of this adaptive process is sketched in Algorithm 1.

\begin{table}[ht!]
\begin{center}
\begin{tabular}{l}
\hline
\rule[0mm]{0mm}{3ex}
\textbf{Algorithm 1: Adaptive procedure} \\    
\hline
\rule[0mm]{0mm}{3ex}
\medskip
\hspace{-0.1cm}\textsc{Step 1}\hspace{0.5cm} Consider the set $X^{0}\equiv X$ of interpolation points\\
\medskip
\textsc{Step 2}\hspace{0.5cm} Fix two positive tolerances (or thresholds) $\theta_{\mbox{refine}}$ and $\theta_{\mbox{coarse}}$, \\
\smallskip
\hskip 1.cm \hskip 1.cm such that $0 < \theta_{\mbox{refine}} < \theta_{\mbox{coarse}}$ \\
\medskip
\textsc{Step 3}\hspace{0.5cm} For $k=0,1,\ldots$ compute the $k$th approximate solution $s_{X^{(k)}}$ \\
\medskip
\hskip 1.cm\textsc{Step 4}\hspace{0.5cm} Define a set $T^{(k)}$ of test points \\
\medskip
\hskip 1.cm\textsc{Step 5}\hspace{0.5cm} Evaluate the residual error $\xi(\boldsymbol{t}_i^{(k)})$ in \eqref{resid}  \\
\medskip
\hskip 1.cm\textsc{Step 6}\hspace{0.5cm} If the error indicator \\
\smallskip
\hskip 1.cm\hspace{1.0cm} \hspace{1.cm} \texttt{i)} $\xi(\boldsymbol{t}_{i}^{(k)}) > \theta_{\mbox{refine}}$, add the test point $\boldsymbol{t}_i^{(k)}$ among the interpolation points \\
\smallskip
\hskip 1.cm\hspace{1.0cm} \hspace{1.cm} \texttt{ii)} $\xi(\boldsymbol{t}_{i}^{(k)}) < \theta_{\mbox{coarse}}$, remove the interpolation node $\boldsymbol{x}_i^{(k)}$ from  $X^{(k)}$ \\
\smallskip
\hskip 1.cm\hspace{1.0cm} \hspace{1.5cm} and put it in the set $X^{(k)}_{\mbox{coarse}}$  \\
\medskip
\hskip 1.cm\hspace{1.5cm} So define the sets\\
\smallskip
\hskip 1.cm\hspace{1.0cm} \hspace{1.cm} $X^{(k)}_{\mbox{refine}} = \{\boldsymbol{t}_i^{(k)}\in T^{(k)} \, : \, \xi(\boldsymbol{t}_{i}^{(k)}) > \theta_{\mbox{refine}}, \, i=1,\ldots,N_{T^{(k)}}\}$\\
\smallskip
\hskip 1.cm\hspace{1.0cm} \hspace{1.cm}	$X^{(k)}_{\mbox{coarse}} = \{\boldsymbol{x}_i^{(k)}\in X^{(k)} \, : \, \xi(\boldsymbol{t}_{i}^{(k)}) < \theta_{\mbox{coarse}}, \, i=1,\ldots,N_{T^{(k)}}\}$\\
\medskip
\hskip 1.cm\hspace{1.5cm} and construct the set \\
\smallskip
\hskip 1.cm\hspace{1.0cm} \hspace{1.cm} $X^{(k+1)} = \left\{ X^{(k)} \cup X^{(k)}_{\mbox{refine}} \right\} \backslash X^{(k)}_{\mbox{coarse}}$\\
\medskip
\hskip 1.cm\textsc{Step 7}\hspace{0.5cm} Stop when $X^{(k)}_{\mbox{refine}} \cup X^{(k)}_{\mbox{coarse}}= \emptyset$\\
\hline
\end{tabular}
\end{center}
\end{table}

\begin{remark}
This adaptive process is based on the computation of the residual \eqref{resid}. It is therefore evident that at each iteration the method requires to create an interpolant and to make some extra evaluations of the target function at the test points. This fact may not be a positive feature when the function evaluation is costly or possibly not available. However, in such a case one might use an alternative approach, which consists in generating for instance a local approximation around the test point, then considering the latter (instead of function value) in \eqref{resid}. Similar strategies have already been studied in e.g. \cite{cav20a,zha17}.
\end{remark}

\subsection{Connection between interpolation and check points} \label{sec:42}

In this subsection we outline the strategy employed for the definition of the sets $X^{(k)}$ and $T^{(k)}$ above. Here we take $k$ fixed, meaning that $k$ is a generic iteration of our adaptive scheme.

In the sequel we describe the connection between the interpolation node set $X^{(k)}$ and the corresponding check point set $T^{(k)}$. By doing that, we focus more in detail on two specific situations that refer to one-dimensional and two-dimensional interpolation.

In 1D case, we start by generating a set $X^{(0)}=X$ of equally spaced points in the domain $\Omega=[a,b]$, $a,b \in \RR$. Then, for $k \geq 0$, we define the set $T^{(k)}$ of test nodes that are the middle points taken from (sorted) interpolation nodes, i.e. $T^{(k)} = \{t_i^{(k)} = 0.5(x_i^{(k)} + x_{i+1}^{(k)}),\ i = 1, \ldots, N^{(k)} - 1\}$.

\begin{figure}[ht!]
\centering
{
\includegraphics[scale=0.54]{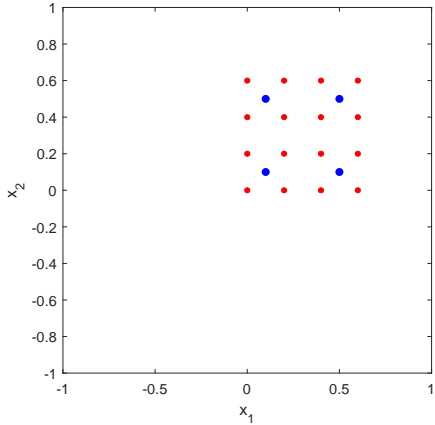}
}
\caption{Example of interpolation nodes (blue) vs check points (red) that refer to a portion of the domain $\Omega=[-1,1]^2$.}
\label{fig:check_vs_nodes}
\end{figure}

In 2D case, we follow the procedure described in \cite{dri07,zha17}. So we start from a set $X^{(0)}=X$ of equally spaced points in the square domain $\Omega=[a,b]^2$, $a,b \in \RR$, and then we update the node set $X^{(k)}$ by applying the adaptive sub-sampling procedure. Now, for $k \geq 0$, we compute the halfway points of $T^{(k)}$ (red dots), as shown in Figure \ref{fig:check_vs_nodes}, where the blue points represent a portion of the set $X^{(k)}$ of interpolation nodes.

\section{Numerical results and discussion} \label{sec:5}

In this section we illustrate the effectiveness of our adaptive algorithms, which are implemented in \textsc{Matlab} for kernel interpolation in one and two dimensions. All programs are run on a laptop with an Intel(R) Core(TM) i7-6500U CPU 2.50 GHz processor with 8GB RAM. 

In the numerical experiments we test our interpolation scheme in order to analyze primarily the behavior of the algorithms in terms of both accuracy and computational efficiency. Then, we also emphasize on some important details for an adaptive method, that is, the number of iterations ($\#$ iter), the final number of points needed for achieving the algorithm convergence ($N_{fin}$) and the conditioning of the interpolation matrix ($\kappa(\mA)$). Furthermore, we compare our results with those of the adaptive algorithms in \cite{dri07}. In our examples, in order to show how the new algorithms work, we consider various types of radial kernels thus involving both infinity and finite regularity RBFs like IMQ, M6, M4 and M2. In our MPLE-RSM we select the shape parameter $\varepsilon$ as discussed in Section \ref{sec:3}. In particular, as suggested in Remark \ref{rem}, the $\varepsilon$-choice via MPLE is determined by the use of \textsc{Matlab} \texttt{fminbnd} minimization. Instead, as regards the comparison of our method with the RSM, we remind that the variable shape parameter selection concerning the RSM derives from the original paper \cite{dri07}.

In order to measure the accuracy of our adaptive method, we compute the $\infty$-norm error or \emph{maximum absolute error} (MAE) given by
\begin{align*}
\mbox{MAE} = ||f-s_{X}||_{\infty}=\max_{1\leq i \leq N_{e}} |f(\boldsymbol{\xi}_i)-s_{X}(\boldsymbol{\xi}_i)|,
\end{align*} 
where the $\boldsymbol{\xi}_i$ forms a suitable set of $N_e$ equally-spaced or gridded evaluation points. Further, by making use of the \textsc{Matlab} \texttt{cond} command we provide an estimate of the condition number \eqref{condnum}, while the efficiency of the adaptive algorithms is assessed by computing the execution (or CPU) time expressed in seconds.

\subsection{Experiments for 1D adaptive interpolation} \label{sec:51}

In this subsection we focus on one-dimensional interpolation. All these tests have been carried by starting from an initial point set $X^{(0)} \equiv X$, which consists of $N^{(0)}=13$ equally-spaced points in the interval $[-1,1]$. The threshold values are usually selected to be $\theta_{\mbox{refine}} = 10^{-5}, 10^{-6}$ and $\theta_{\mbox{coarse}}=10^{-8}, 10^{-9}$. However, in our comparison between the residual sub-sampling algorithms, the refinement threshold $\theta_{\mbox{refine}}$ is often modified by keeping $\theta_{\mbox{coarse}}$ fixed.

In order to validate in depth our adaptive algorithms, we consider the following three benchmark target (or test) functions: 
\begin{align*}
f_1(x) = \frac{1}{1+25x^2}, \qquad  f_2(x) = \tanh(60x-0.01), \qquad f_3(x) = \frac{3}{8}\left[\cos\left(\left(x+1\right)^2-3\right)\right]^4,
\end{align*}
where $f_1$ is the well-known Runge function, $f_2$ represents the hyperbolic tan function, and $f_3$ denotes a univariate restriction of the so-called valley function (see \cite{cav21b,dri07,zha17}).

\begin{figure}[ht!]
\centering
{
\includegraphics[scale=0.54]{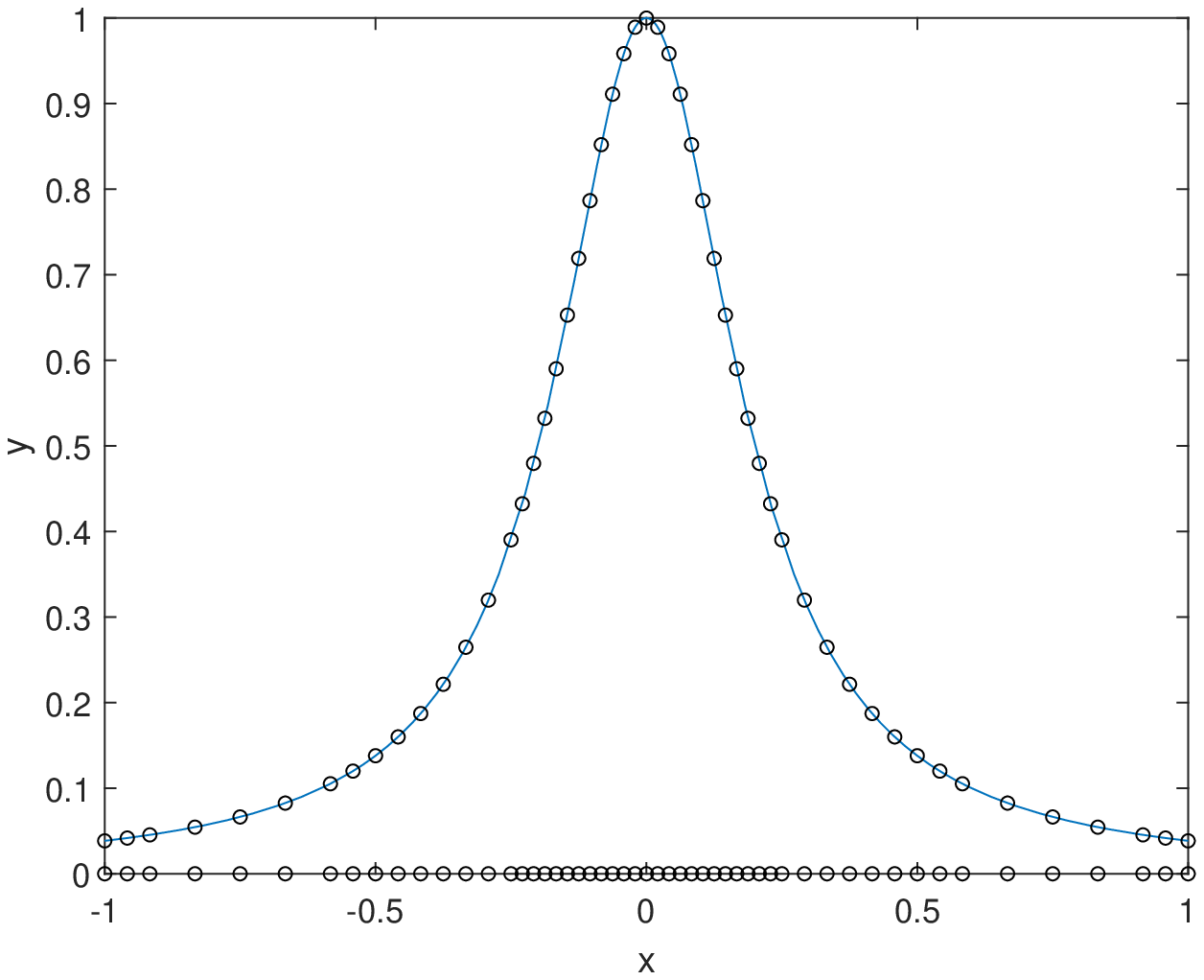}
\includegraphics[scale=0.54]{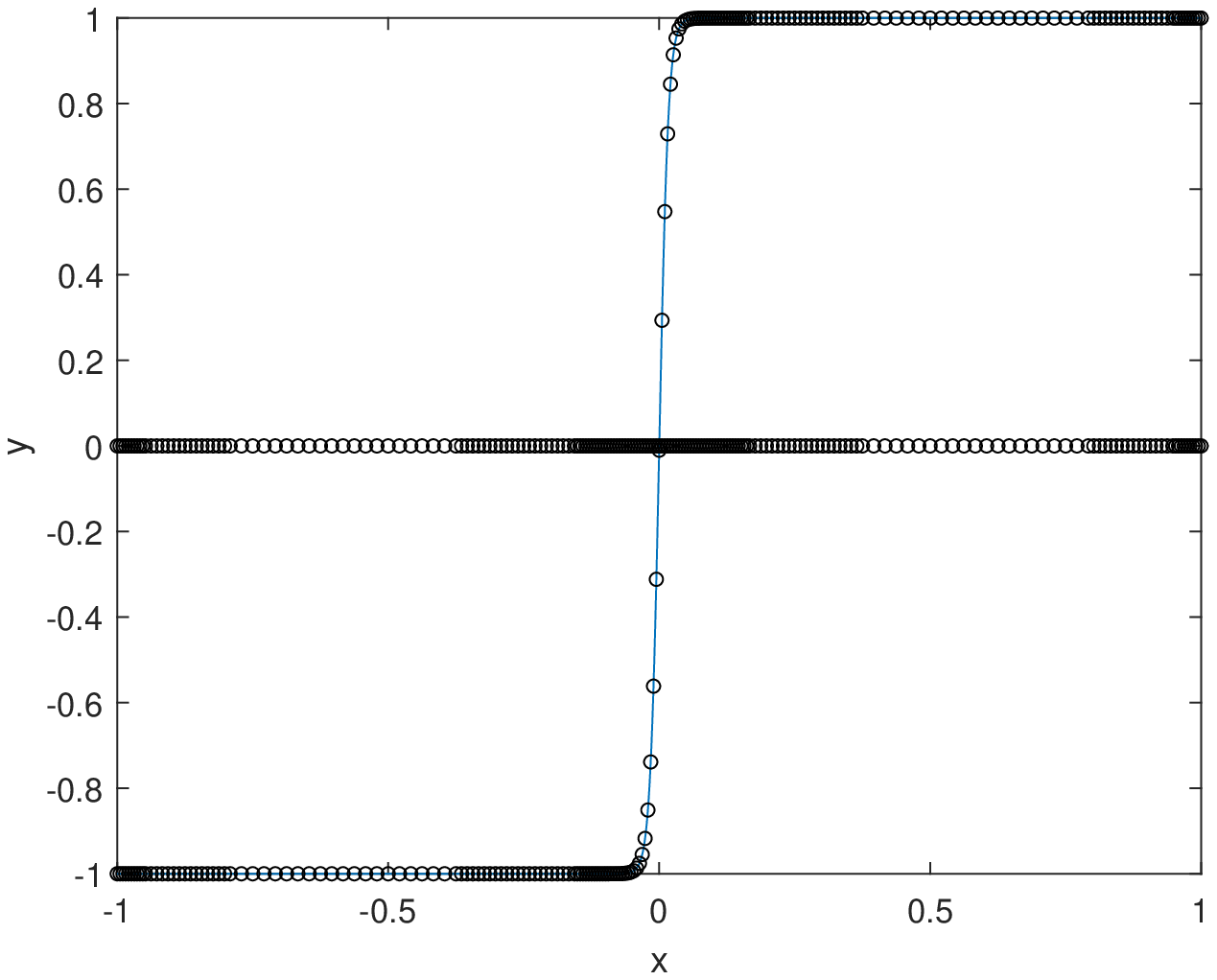}\\
\includegraphics[scale=0.54]{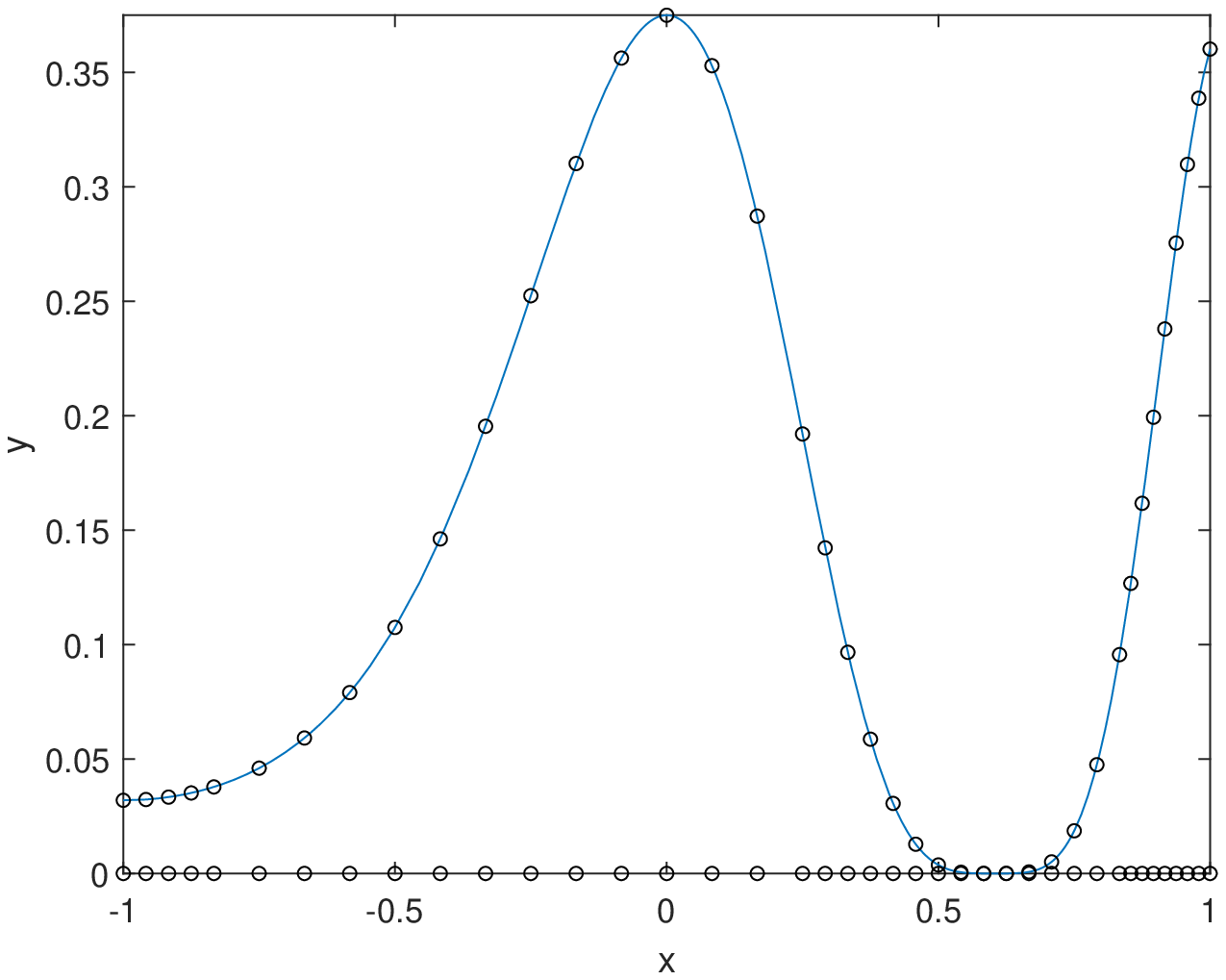}
}
\caption{Final point distribution using MPLE-RSM. 1D adaptive interpolation for $f_1$ with M4 (top-left), $f_2$ with IMQ (top-right) and $f_3$ with M6 (bottom). The chosen parameters are given in Tables \ref{tab:1}, \ref{tab:2} and \ref{tab:3}, respectively.}
\label{fig:test_1d}
\end{figure}

In Figure \ref{fig:test_1d} we show some final point distributions obtained by adaptive interpolation. These results are three examples of possible application of MPLE-RSM, which considers different choices of kernel for $f_1$, $f_2$ and $f_3$. For the Runge function $f_1$ (top-left) we observe that the points cluster close to the boundaries where approximation turns out to be more challenging due to the one-side nature of the information, and at the origin in which the target function changes more rapidly. In case of hyperbolic tan function (top-right) we observe as the points distribute around the steepest part of $f_2$. A similar behavior also occurs for $f_3$, since the nodes tends to gather near the boundaries and the steepest areas of this function (bottom).

In Tables \ref{tab:1}, \ref{tab:2} and \ref{tab:3}, we present the results obtained by applying the adaptive MPLE-RSM. From these tables we have some information about the number of iterations and the final number of points required for convergence. In particular, we note that the MAE is always very close to the prescribed value of $\theta_{\mbox{refine}}$. This fact suggests that our method is working well. Moreover, the automatic selection of the optimal shape parameter also permits to control the condition number (in all cases $< 10^{+15}$, even for $C^{\infty}$ kernels), at the same time guaranteeing a high level of accuracy of the numerical method. Indeed, in the various situations, conditioning and CPU times assume quite similar values. It is, however, important to notice that the steep variation in the middle of $f_2$ leads to a remarkable increase in the final number of points, which the adaptive process needs to achieve the $\theta_{\mbox{refine}}$ threshold. In the latter case, in fact, for any choice of the kernel $N_{fin}$ assumes a value larger than $100$, while for $f_1$ and $f_3$ the value of $N_{fin}$ is always smaller than $100$.

\begin{table}[th!]
\begin{center}
{ 
\begin{tabular}{|c||c|c|c|c|c|} \hline
				\rule[-2mm]{0mm}{7mm}
kernel &$\#$ iter&$N_{fin}$&MAE&$\kappa(\mA)$&time\\ \hline
			  \rule[0mm]{0mm}{3ex}
IMQ &4 &51&2.1e-07&3.0e+12&0.6\\
			  \rule[0mm]{0mm}{3ex}
M6  &10&50&6.7e-07&1.5e+11&0.8\\
			  \rule[0mm]{0mm}{3ex}
M4  &12&54&8.9e-07&8.4e+09&1.1\\
			  \rule[0mm]{0mm}{3ex}
M2  &6 &99&9.3e-07&4.4e+08&0.7\\
\hline
\end{tabular}
}
\caption{Results obtained by applying the MPLE-RSM with $\theta_{\mbox{refine}} = 10^{-6}$ and $\theta_{\mbox{coarse}} = 10^{-8}$ for $f_1$.}
\label{tab:1}
\end{center}
\end{table}

\begin{table}[th!]
\begin{center}
{ 
\begin{tabular}{|c||c|c|c|c|c|} \hline
				\rule[-2mm]{0mm}{7mm}
kernel &$\#$ iter&$N_{fin}$&MAE&$\kappa(\mA)$&time\\ \hline
					  \rule[0mm]{0mm}{3ex}
IMQ &6  &194 &3.5e-06 &5.1e+13 &0.7 \\
					  \rule[0mm]{0mm}{3ex}
M6  &14 &108 &8.9e-06 &1.0e+11 &1.7 \\
					  \rule[0mm]{0mm}{3ex}
M4  &11 &146 &9.8e-06 &1.0e+09 &1.0 \\
					  \rule[0mm]{0mm}{3ex}
M2  &8  &154 &9.8e-06 &1.9e+08 &0.7 \\
\hline
\end{tabular}
}
\caption{Results obtained by applying the MPLE-RSM with $\theta_{\mbox{refine}} = 10^{-5}$ and $\theta_{\mbox{coarse}} = 10^{-8}$ for $f_2$.}
\label{tab:2}
\end{center}
\end{table}

\begin{table}[th!]
\begin{center}
{ 
\begin{tabular}{|c||c|c|c|c|c|} \hline
				\rule[-2mm]{0mm}{7mm}
kernel &$\#$ iter&$N_{fin}$&MAE&$\kappa(\mA)$&time\\ \hline
					  \rule[0mm]{0mm}{3ex}
IMQ &4 &29 &2.3e-06 &2.3e+14 &0.6 \\
					  \rule[0mm]{0mm}{3ex}
M6  &4 &40 &5.6e-06 &6.9e+11 &0.4 \\
					  \rule[0mm]{0mm}{3ex}
M4  &5 &43 &1.0e-05 &1.8e+10 &0.6 \\
					  \rule[0mm]{0mm}{3ex}
M2  &7 &56 &9.2e-06 &1.1e+10 &0.6 \\
\hline
\end{tabular}
}
\caption{Results obtained by applying the MPLE-RSM with $\theta_{\mbox{refine}} = 10^{-5}$ and $\theta_{\mbox{coarse}} = 10^{-8}$ for $f_3$.}
\label{tab:3}
\end{center}
\end{table}

Finally, in Tables \ref{tab:4}, \ref{tab:5} and \ref{tab:6} we compare our adaptive MPLE-RSM with the RSM proposed in \cite{dri07}. Thus, we report the results obtained by varying the refinement threshold $\theta_{\mbox{refine}}$ for the M6 kernel. This analysis enables us to make some general observations. When the value of $\theta_{\mbox{refine}}$ is \lq\lq small\rq\rq\ the MPLE-RSM achieves convergence much faster than RSM. Furthermore, while for \lq\lq large\rq\rq\ values of $\theta_{\mbox{refine}}$ the number of points is pretty similar for the two methods, the MPLE-RSM usually needs much less points than RSM. These remarks are true for each of the target functions. However, from these experiments we can note that for the RSM \cite{dri07} is not always possible to get any result, and so in the tables we denote this issue with the symbol --. Another drawback of the RSM is then due to severe difficulties in tuning the variable shape parameters. Unlike the MPLE-RSM, where the shape parameter choice is automatic for any radial kernel, the RSM needs user's action case-by-case thus making the $\varepsilon$-selections quite hard. This fact is particularly evident, either when the target function (and, as a consequence, the approximation problem) is quite complex, or the interpolation problem requires to be solved by kernels that have different degrees of smoothness.

\begin{table}[th!]
\begin{center}
{ 
\begin{tabular}{|c||c|c||c|c|} \hline
			  \rule[-2mm]{0mm}{7mm}
\multirow{3}{*}{$\theta_{\mbox{refine}}$} & \multicolumn{2}{|c||}{RSM} & \multicolumn{2}{|c|}{MPLE-RSM}\\
\cline{2-5}
				\rule[-2mm]{0mm}{7mm}
    &$N_{fin}$&time&$N_{fin}$&time\\ \hline
					  \rule[0mm]{0mm}{3ex}
1e-04 &25&0.2&25&0.2\\
			  \rule[0mm]{0mm}{3ex}
1e-05 &97&0.6&35&0.2\\
			  \rule[0mm]{0mm}{3ex}
1e-06 &224&3.1&49&0.3\\
			  \rule[0mm]{0mm}{3ex}
1e-07 &--&--&67&0.4\\
\hline
\end{tabular}
}
\caption{Comparison between RSM \cite{dri07} and MPLE-RSM obtained by using the M6 kernel with $\theta_{\mbox{coarse}} = 10^{-9}$ for $f_1$.}
\label{tab:4}
\end{center}
\end{table}

\begin{table}[th!]
\begin{center}
{ 
\begin{tabular}{|c||c|c||c|c|} \hline
\rule[-2mm]{0mm}{7mm}
\multirow{3}{*}{$\theta_{\mbox{refine}}$} & \multicolumn{2}{|c||}{RSM} & \multicolumn{2}{|c|}{MPLE-RSM}\\
\cline{2-5}
				\rule[-2mm]{0mm}{7mm}
    &$N_{fin}$&time&$N_{fin}$&time\\ \hline
							  \rule[0mm]{0mm}{3ex}
1e-03 &73 &0.7&82 &0.5\\
					  \rule[0mm]{0mm}{3ex}
1e-04 &477&3.0&112&0.8\\
					  \rule[0mm]{0mm}{3ex}
1e-05 &407&4.3&108&1.7\\
					  \rule[0mm]{0mm}{3ex}
1e-06 &--&--&129&4.5\\
\hline
\end{tabular}
}
\caption{Comparison between RSM \cite{dri07} and MPLE-RSM obtained by using the M6 kernel with $\theta_{\mbox{coarse}} = 10^{-8}$ for $f_2$.}
\label{tab:5}
\end{center}
\end{table}

\begin{table}[th!]
\begin{center}
{ 
\begin{tabular}{|c||c|c||c|c|} \hline
\rule[-2mm]{0mm}{7mm}
\multirow{3}{*}{$\theta_{\mbox{refine}}$} & \multicolumn{2}{|c||}{RSM} & \multicolumn{2}{|c|}{MPLE-RSM}\\
\cline{2-5}
				\rule[-2mm]{0mm}{7mm}
    &$N_{fin}$&time&$N_{fin}$&time\\ \hline
							  \rule[0mm]{0mm}{3ex}
1e-03 &20 &0.3&20&0.3\\
					  \rule[0mm]{0mm}{3ex}
1e-04 &92 &0.7&30&0.6\\
					  \rule[0mm]{0mm}{3ex}
1e-05 &184&1.9&40&0.4\\
					  \rule[0mm]{0mm}{3ex}
1e-06 &-- &-- &40&0.7\\
\hline
\end{tabular}
}
\caption{Comparison between RSM \cite{dri07} and MPLE-RSM obtained by using the M6 kernel with $\theta_{\mbox{coarse}} = 10^{-8}$ for $f_3$.}
\label{tab:6}
\end{center}
\end{table}

\subsection{Experiments for 2D adaptive interpolation} \label{sec:52}

In this subsection we consider the two-dimensional adaptive interpolation algorithm. These experiments have been run by taking an initial point set $X^{(0)} \equiv X$, containing $N^{(0)}=320$ uniformly distributed points on $[-1,1]^2$. In order to test the node refinement process, as refinement thresholds we choose some values of $\theta_{\mbox{refine}} \in [10^{-6}, 10^{-3}]$, while the coarsening tolerance is assumed to be fixed, i.e., $\theta_{\mbox{coarse}}=10^{-8}$. Furthermore, as in the 1D case above, we compare the numerical results obtained by using the classical RSM in \cite{dri07} with the new MPLE-RSM. In doing that, we analyze the behavior of the two algorithms by varying the refinement threshold.

In our tests we analyze the performance of our algorithms taking the data values by three test functions. The former is known as a Franke-type function \cite{zha17}, and its analytic expression is
\begin{align*}
   f_4(x,y) &= \exp\left[-0.1\left(x^2 + y^2\right)\right] + \exp\left[-5\left((x - 0.5)^2 + (y - 0.5)^2\right)\right] \\
       &+ \exp\left[-15((x+0.2)^2 + (y+0.4)^2)\right] + \exp\left[-9\left((x+0.8)^2 + (y-0.8)^2\right)\right].
\end{align*}
The latter is a hyperbolic tan function \cite{dri07} of the form
\begin{align*}
   f_5(x,y) = -0.4\tanh(20xy)+0.6,
\end{align*}
while the last one is an exponential function \cite{zha17} given by
\begin{align*}
   f_6(x,y) = \exp\left[-60((x-0.35)^2+(y-0.25)^2)\right]+0.2.
\end{align*}

\begin{figure}[ht!]
\centering
{
\includegraphics[scale=0.54]{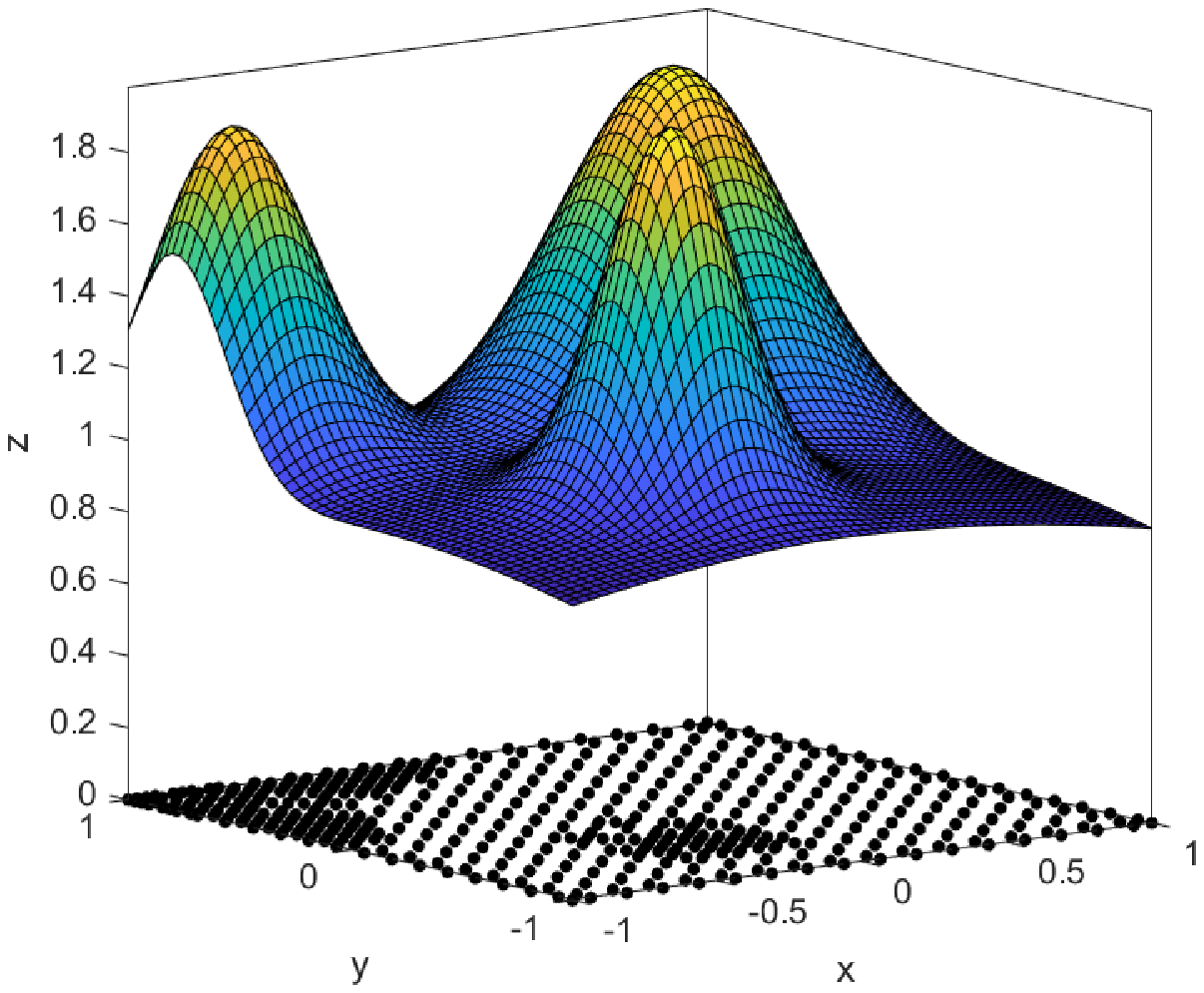}
\includegraphics[scale=0.54]{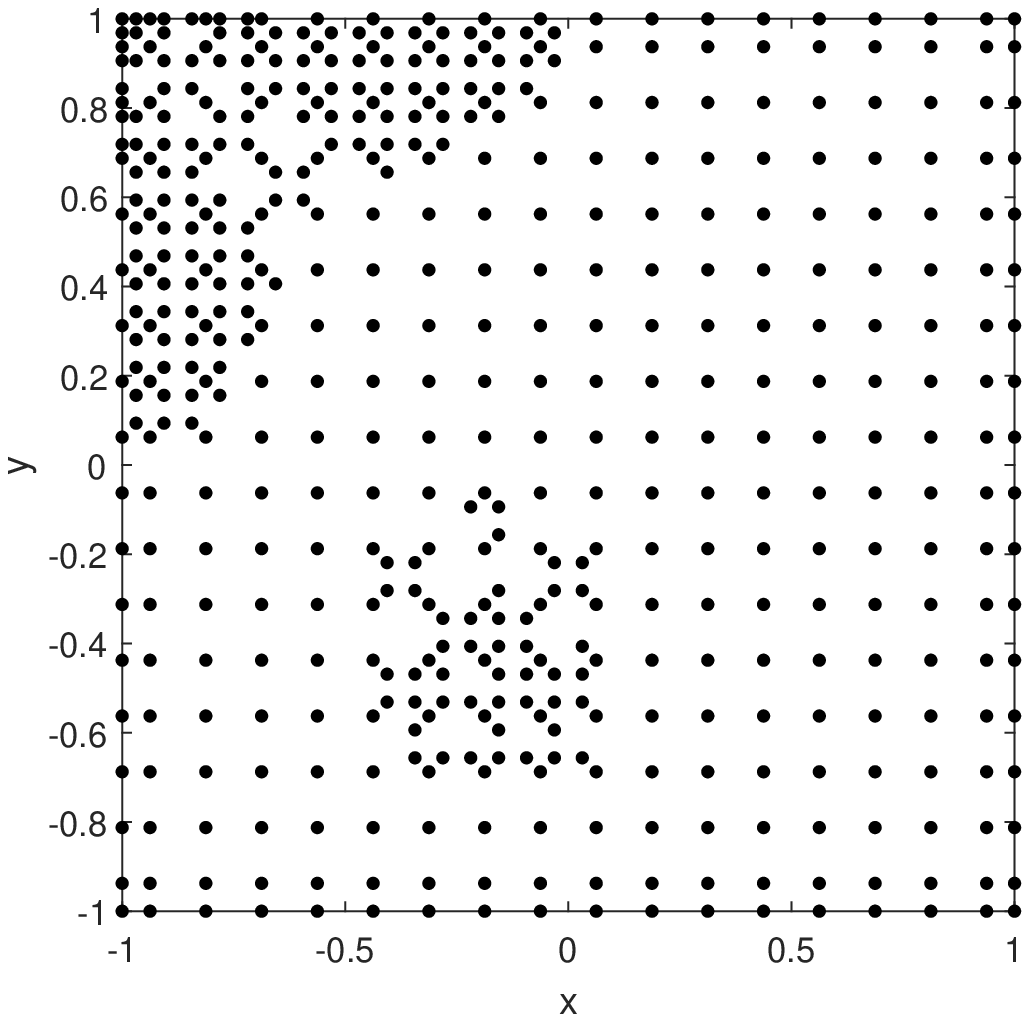}\\
\includegraphics[scale=0.54]{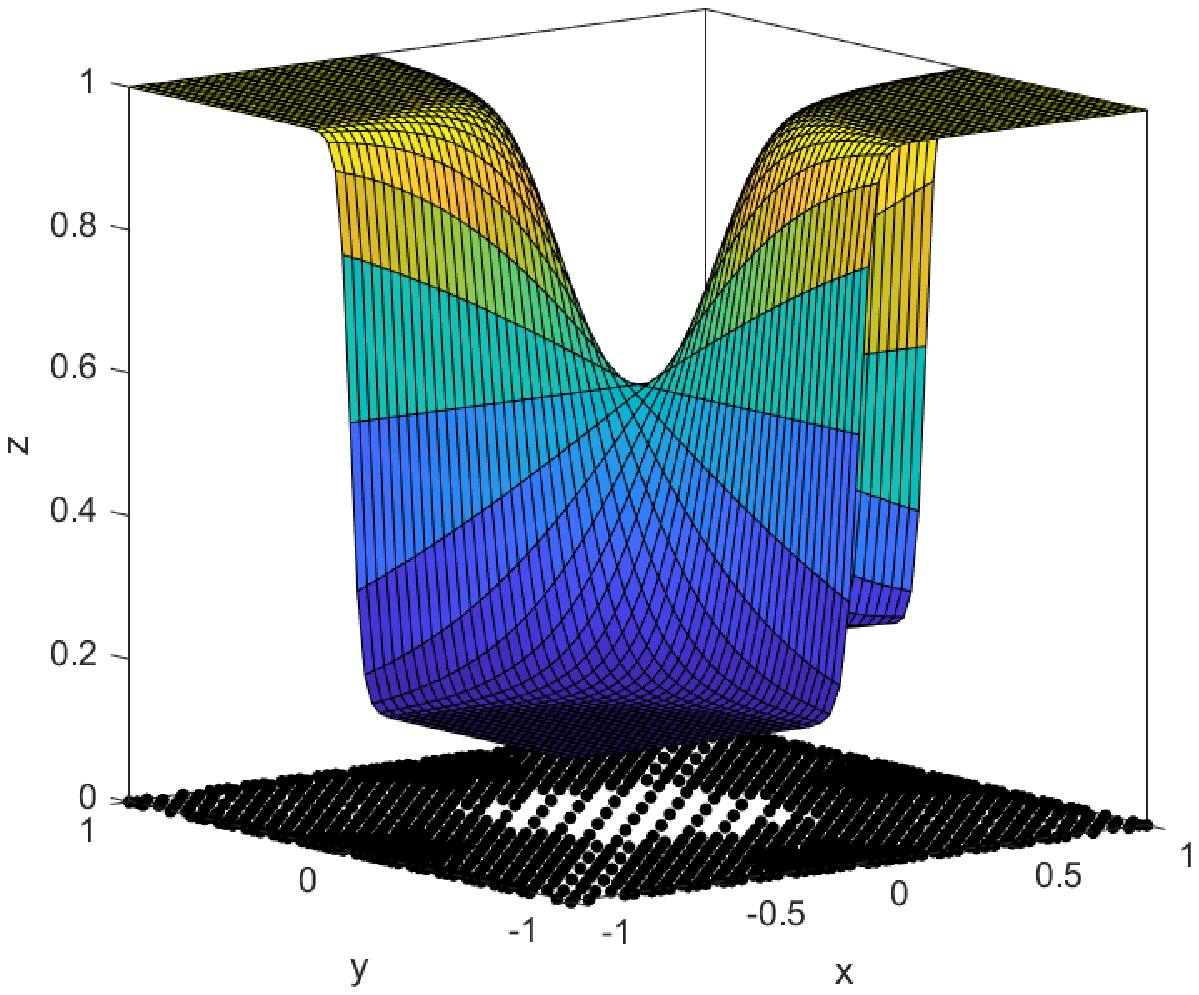}
\includegraphics[scale=0.54]{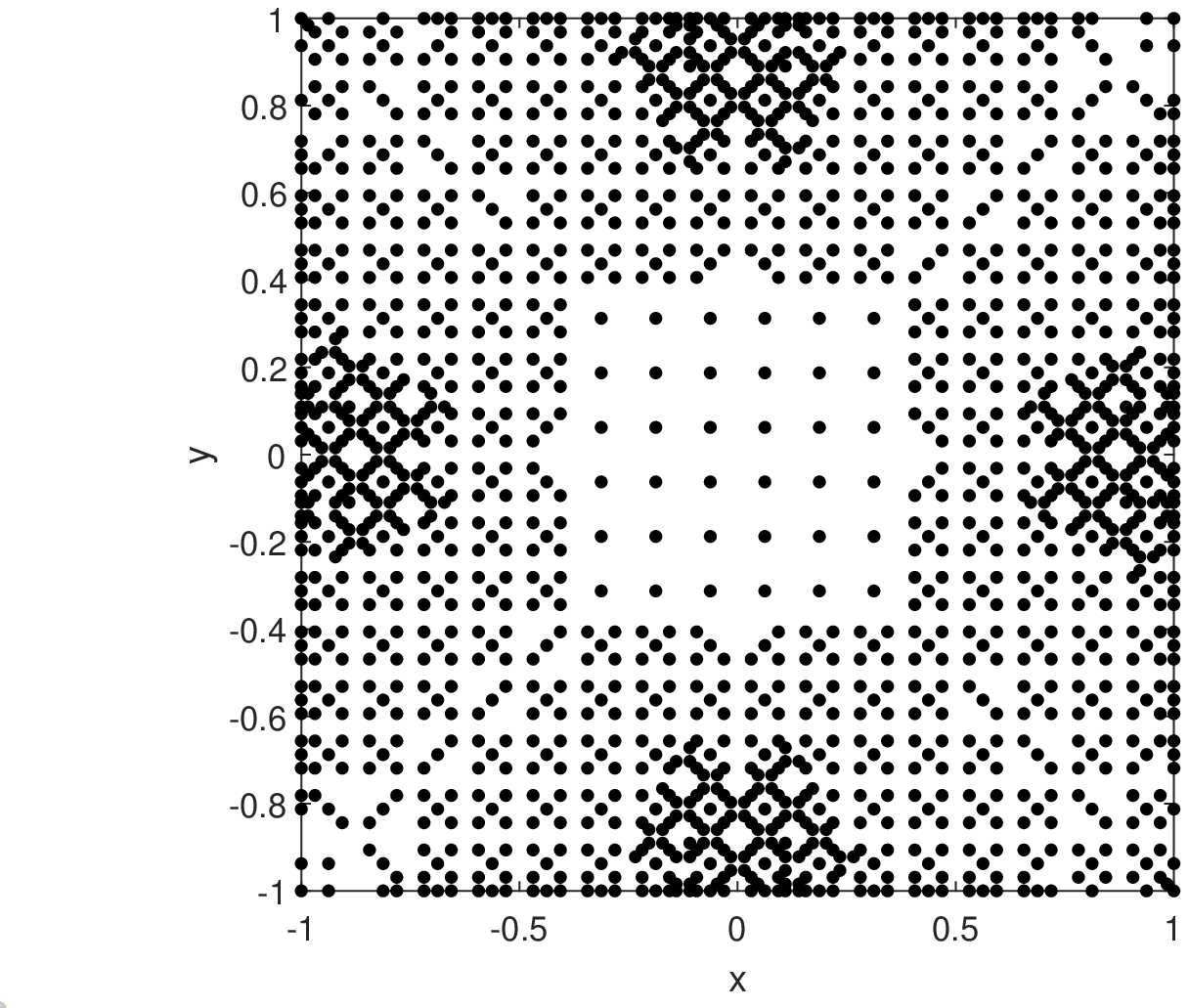}\\
\includegraphics[scale=0.54]{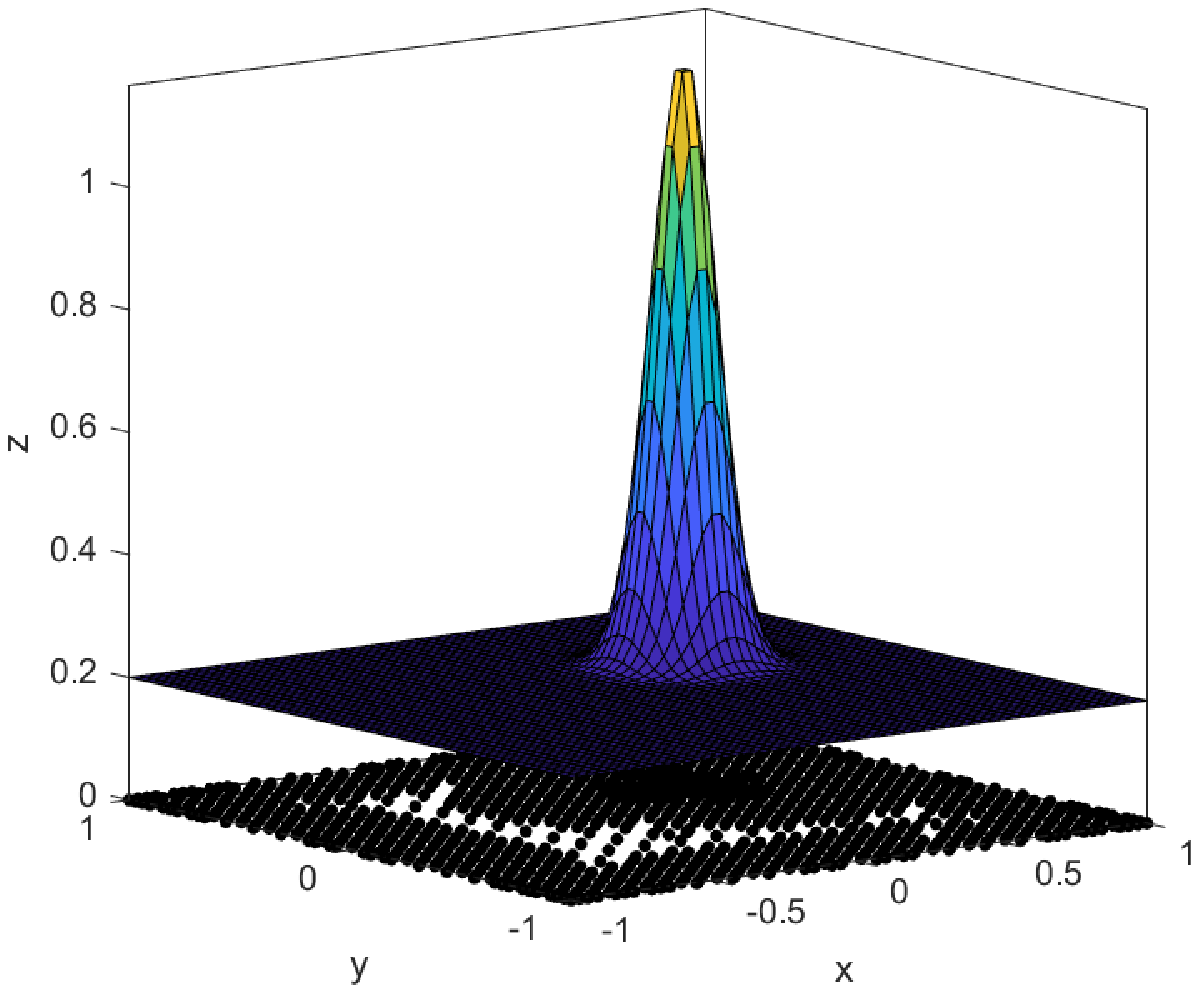}
\includegraphics[scale=0.54]{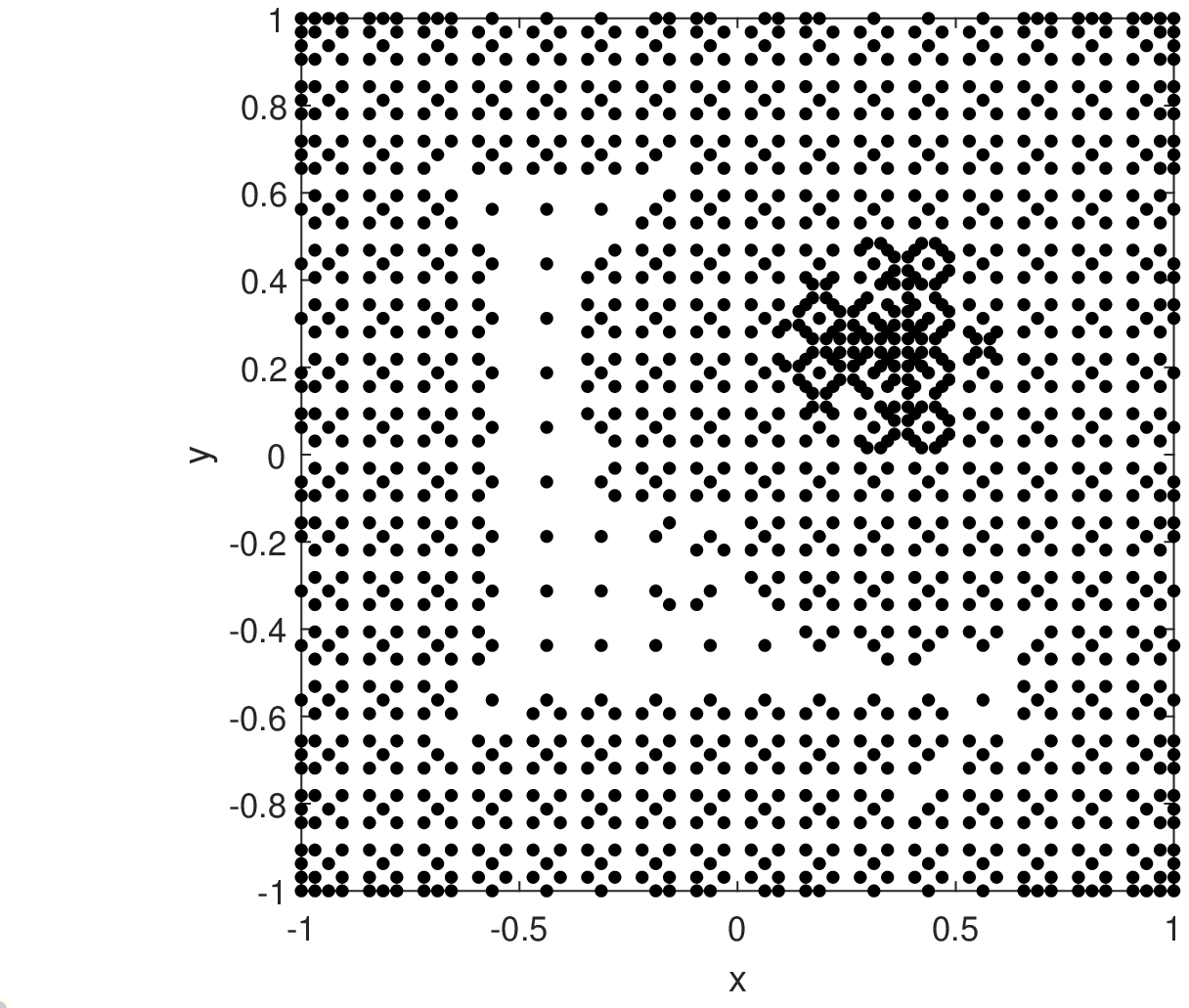}
}
\caption{Final point distribution using MPLE-RSM. 2D adaptive interpolation for $f_4$ with M4 (top), $f_5$ with IMQ (center) and $f_6$ with M6 (bottom). The chosen parameters are given in Tables \ref{tab:7}, \ref{tab:8} and \ref{tab:9}, respectively.}
\label{fig:test_2d}
\end{figure}

In Figure \ref{fig:test_2d} we give some graphical representations of the final point distribution, which are obtained by applying our adaptive MPLE-RSM algorithm. On the left, we report the graphs of the interpolating functions for $f_4$ (top), $f_5$ (center) and $f_6$ (bottom), also showing on the $xy$-plane at level $z=0$ the final nodes deriving from the adaptive interpolation. On the right, instead, we explicitly show the 2D view of the point distributions. To give a generic idea about the flexibility of this kernel based interpolation approach, we depict the results with radial kernels of various regularity; as an example, in our work, we report the results attained in the following cases: M4 for $f_4$, IMQ for $f_5$ and M6 for $f_6$. More specifically, we can observe that for the Franke-type function $f_4$ the MPLE-RSM locates points in regions of rapid variation. The adaptive algorithm behaves in a similar way also when we consider the hyperbolic tan function $f_5$ and the exponential one $f_6$. In both cases the residual sub-sampling scheme puts more points in the domain where the functions change quickly or are picked.

In Tables \ref{tab:7}, \ref{tab:8} and \ref{tab:9}, we illustrate the numerical results to see how the adaptive MPLE-RSM works when it is applied to solve some relevant approximation problems. As already shown in 1D interpolation above, in the tables we provide a detailed summary concerning  the execution of the adaptive algorithm. In particular, we report the number of iterations needed to achieve convergence and the corresponding final number of points. Also in these 2D tests, we can observe that the MAE assumes values that are close -- and however always lower -- to the prefixed refinement threshold $\theta_{\mbox{refine}}$. Such data point out that the MPLE-RSM is able to make accurate predictions. In addition, we note that the near-optimal determination of the shape parameter via the MPLE technique enables us to control the kernel matrix conditioning. Indeed, in all reported results the condition number is smaller than $10^{+16}$, even for $C^{\infty}$ RBFs like IMQ. As regards the CPU times we highlight that in this adaptive context the selected kernel can influence the execution time of the MPLE-RSM, which is subjected to an automatic addition or removal of points. Moreover, although the quality of results depends on the choice of the refinement threshold $\theta_{\mbox{refine}}$ and the complexity of the target function to be approximated, it is important to note as the number of refinement nodes remains relatively low. This fact is better explained in the following, where we compare the MPLE-RSM with the original RSM \cite{dri07}.

\begin{table}[th!]
\begin{center}
{ 
\begin{tabular}{|c||c|c|c|c|c|} \hline
				\rule[-2mm]{0mm}{7mm}
kernel &    $\#$ iter&$N_{fin}$&MAE&$\kappa(\mA)$&time\\ \hline
							  \rule[0mm]{0mm}{3ex}
IMQ &1 &318 &5.3e-05 &9.6e+13 &0.6 \\
					  \rule[0mm]{0mm}{3ex}
M6  &3 &398 &5.7e-05 &1.6e+12 &1.9 \\
					  \rule[0mm]{0mm}{3ex}
M4  &3 &491 &7.1e-05 &1.6e+11 &0.8 \\
\hline
\end{tabular}
}
\caption{Results obtained by applying the MPLE-RSM with $\theta_{\mbox{refine}} = 10^{-4}$ and $\theta_{\mbox{coarse}} = 10^{-8}$ for $f_4$.}
\label{tab:7}
\end{center}
\end{table}

\begin{table}[th!]
\begin{center}
{ 
\begin{tabular}{|c||c|c|c|c|c|} \hline
					  \rule[-2mm]{0mm}{7mm}
kernel  &$\#$ iter&$N_{fin}$&MAE&$\kappa(\mA)$&time\\ \hline
							  \rule[0mm]{0mm}{3ex}
IMQ &4 &1522 &5.8e-04 &6.0e+10 &7.6 \\
			  \rule[0mm]{0mm}{3ex}
M6  &3 &1442 &7.2e-04 &8.1e+08 &9.9 \\
			  \rule[0mm]{0mm}{3ex}
M4  &3 &1300 &8.4e-04 &2.7e+08 &4.0 \\
\hline
\end{tabular}
}
\caption{Results obtained by applying the MPLE-RSM with $\theta_{\mbox{refine}} = 10^{-3}$ and $\theta_{\mbox{coarse}} = 10^{-8}$ for $f_5$.}
\label{tab:8}
\end{center}
\end{table}

\begin{table}[th!]
\begin{center}
{ 
\begin{tabular}{|c||c|c|c|c|c|} \hline
				\rule[-2mm]{0mm}{7mm}
kernel  &$\#$ iter&$N_{fin}$&MAE&$\kappa(\mA)$&time\\ \hline
			  \rule[0mm]{0mm}{3ex}
IMQ &2 &1259 &3.4e-06 &5.7e+15 &3.1 \\
			  \rule[0mm]{0mm}{3ex}
M6  &4 &1368 &6.1e-06 &2.3e+12 &18.1 \\
			  \rule[0mm]{0mm}{3ex}
M4  &4 &1428 &9.3e-06 &9.3e+10 &8.6 \\
\hline
\end{tabular}
}
\caption{Results obtained by applying the MPLE-RSM with $\theta_{\mbox{refine}} = 10^{-5}$ and $\theta_{\mbox{coarse}} = 10^{-8}$ for $f_6$.}
\label{tab:9}
\end{center}
\end{table}

Therefore, in order to point out the benefit deriving from the use of the new adaptive algorithm, we conclude this numerical section by making a comparison between MPLE-RSM and RSM \cite{dri07}. In Tables \ref{tab:10}, \ref{tab:11} and \ref{tab:12} we report the results obtained for various choices of the refinement parameter $\theta_{\mbox{refine}}$. In this study we also diversify the type of radial kernel and in particular, for the sake of brevity, we show the algorithm behavior by employing the IMQ for $f_4$, the M2 for $f_5$ and the M6 for $f_6$. From these tables, in which we report the final number of interpolation nodes required to meet both tolerances and the CPU times, it turns out to be undeniable the improvement provided by the novel adaptive MPLE-RSM. As a matter of fact, for any value of $\theta_{\mbox{refine}}$, the MPLE-RSM is much faster than the RSM. This good result in terms of computational efficiency is due to a better ability of our new approach in the selection of the shape parameter via a MPLE based strategy. Consequently, the numerical method achieves convergence earlier and, at the same time, a smaller number of points is required. Furthermore, we can also note a growing enhancement in the performance of the MPLE-RSM (compared to the RSM), when the value of $\theta_{\mbox{refine}}$ becomes more demanding, i.e., the latter is assumed to be smaller and smaller. In some specific cases, then, it is even possible that the RSM does not converge (see the last row of Table \ref{tab:10}). In conclusion, from these numerical experiments for 2D adaptive interpolation, we have once more put in evidence a concrete difficulty in the initialization of shape parameters for the RSM. From several tests it is evident how this (nonautomatic) choice can modify the performance of the algorithm, but the need of a preliminary user action makes the entire scheme highly unstable. However, to perform the comparisons discussed in this work, for the RSM we set as starting values of the shape parameter $\varepsilon = 3$ in Tables \ref{tab:10} and \ref{tab:11}, and $\varepsilon = 4$ in Table \ref{tab:12}.

\begin{table}[th!]
\begin{center}
{ 
\begin{tabular}{|c||c|c|c|c|} \hline
\rule[-2mm]{0mm}{7mm}
\multirow{3}{*}{$\theta_{\mbox{refine}}$} & \multicolumn{2}{|c|}{RSM} & \multicolumn{2}{|c|}{MPLE-RSM}\\
\cline{2-5}
\rule[-2mm]{0mm}{7mm}
    &$N_{fin}$&time&$N_{fin}$&time\\ \hline
					  \rule[0mm]{0mm}{3ex}
1e-04 &612 &1.4 &318 &0.4\\
					  \rule[0mm]{0mm}{3ex}
5e-05 &927 &2.3 &321 &0.8\\
					  \rule[0mm]{0mm}{3ex}
1e-05 &2430 &30.2 &328 &1.3\\
					  \rule[0mm]{0mm}{3ex}
5e-06 &3850 &146.4 &286 &2.3\\
					  \rule[0mm]{0mm}{3ex}
1e-06 &-- &-- &303 &5.0\\
\hline
\end{tabular}
}
\caption{Comparison between RSM \cite{dri07} and MPLE-RSM obtained by using the IMQ kernel with $\theta_{\mbox{coarse}} = 10^{-8}$ for $f_4$.}
\label{tab:10}
\end{center}
\end{table}

\begin{table}[th!]
\begin{center}
{ 
\begin{tabular}{|c||c|c|c|c|} \hline
\rule[-2mm]{0mm}{7mm}
\multirow{3}{*}{$\theta_{\mbox{refine}}$} & \multicolumn{2}{|c|}{RSM} & \multicolumn{2}{|c|}{MPLE-RSM}\\
\cline{2-5}
\rule[-2mm]{0mm}{7mm}
    &$N_{fin}$&time&$N_{fin}$&time\\ \hline
					  \rule[0mm]{0mm}{3ex}
1e-03 &1280 &7.4 &950 &3.0\\
					  \rule[0mm]{0mm}{3ex}
8e-04 &1476 &6.6 &1038 &4.3\\
					  \rule[0mm]{0mm}{3ex}
6e-04 &1750 &9.4 &1172 &4.4\\
					  \rule[0mm]{0mm}{3ex}
4e-04 &2120 &14.8 &1390 &8.9\\
					  \rule[0mm]{0mm}{3ex}
2e-04 &3178 &117.7 &1952 &19.5\\
					  \rule[0mm]{0mm}{3ex}
1e-04 &4732 &289.0 &2406 &41.1\\
\hline
\end{tabular}
}
\caption{Comparison between RSM \cite{dri07} and MPLE-RSM obtained by using the M2 kernel with $\theta_{\mbox{coarse}} = 10^{-8}$ for $f_5$.}
\label{tab:11}
\end{center}
\end{table}

\begin{table}[th!]
\begin{center}
{ 
\begin{tabular}{|c||c|c|c|c|} \hline
\rule[-2mm]{0mm}{7mm}
\multirow{3}{*}{$\theta_{\mbox{refine}}$} & \multicolumn{2}{|c|}{RSM} & \multicolumn{2}{|c|}{MPLE-RSM}\\
\cline{2-5}
\rule[-2mm]{0mm}{7mm}
    &$N_{fin}$&time&$N_{fin}$&time\\ \hline
					  \rule[0mm]{0mm}{3ex}
1e-03 &566 &3.0 &550 &1.4\\
					  \rule[0mm]{0mm}{3ex}
5e-04 &643 &4.3 &658 &1.8\\
					  \rule[0mm]{0mm}{3ex}
1e-04 &978 &7.5 &1215&4.6\\
					  \rule[0mm]{0mm}{3ex}
5e-05 &1186&16.2 &1272&10.5\\
					  \rule[0mm]{0mm}{3ex}
1e-05 &2933&65.1&1368&18.1\\
\hline
\end{tabular}
}
\caption{Comparison between RSM \cite{dri07} and MPLE-RSM obtained by using the M6 kernel with $\theta_{\mbox{coarse}} = 10^{-8}$ for $f_6$.}
\label{tab:12}
\end{center}
\end{table}



\section{Conclusions and future work} \label{sec:6}

In this paper we solved two open problems in \cite{dri07}. In fact, though the original method provides an effective adaptive scheme, it does not guarantee the invertibility of the interpolation matrix. This issue is essentially due to the variable shape parameter selection in the RSM, since the latter is characterized by a different choice of the shape parameter at every node. On the contrary, the new approach based on the MPLE criterion enables us on the one hand to make an optimal choice of the shape parameter associated with the kernel, and on the other one to guarantee existence and uniqueness of the RBF interpolation. Furthermore, the MPLE-RSM is thus automatically applicable to any kind of kernel, while the basic RSM needs a quite hard action of the user for the shape parameter selection. 

As a future work we aim to improve the method proposed in this work, also for the solution of boundary value problems and partial differential equations.


\section*{Acknowledgments}
This work was partially supported by the INdAM-GNCS 2020 research project \lq\lq Multivariate approximation and functional equations for numerical modeling\rq\rq\ and by the 2020 projects \lq\lq Models and numerical methods in approximation, in applied sciences and in life sciences\rq\rq and \lq\lq Mathematical methods in computational sciences\rq\rq\ funded by the Department of Mathematics \lq\lq Giuseppe Peano\rq\rq\ of the University of Turin. This research has been accomplished within the RITA \lq\lq Research ITalian network on Approximation\rq\rq\ and the UMI Group TAA \lq\lq Approximation Theory and Applications\rq\rq. 







\begin{thebibliography}{99}








\bibitem{beh02} J. Behrens, A. Iske, Grid-free adaptive semi-Lagrangian advection using radial basis functions, Comput. Math. Appl. 43 (2002) 319--327.

\bibitem{boz02} M. Bozzini, L. Lenarduzzi, R. Schaback, Adaptive interpolation by scaled multiquadrics, Adv. Comput. Math. 16 (2002) 375--387.

\bibitem{buh03} M.D. Buhmann, Radial Basis Functions: Theory and Implementation, Cambridge Monogr. Appl. Comput. Math., vol. 12, Cambridge Univ. Press, Cambridge, 2003.












\bibitem{cav20a} R. Cavoretto, A. De Rossi, A two-stage adaptive scheme based on RBF collocation for solving elliptic PDEs, Comput. Math. Appl. 79 (2020) 3206--3222.

\bibitem{cav20b} R. Cavoretto, A. De Rossi, An adaptive LOOCV-based refinement scheme for RBF collocation methods over irregular domains, Appl. Math. Lett. 103 (2020) 106178.

\bibitem{cav20c} R. Cavoretto, A. De Rossi, Error indicators and refinement strategies for solving Poisson problems through a RBF partition
of unity collocation scheme, Appl. Math. Comput. 369 (2020) 124824.


\bibitem{cav21b} R. Cavoretto, Adaptive radial basis function partition of unity interpolation: A bivariate algorithm for unstructured data, J. Sci. Comput. 87 (2021) 41.

\bibitem{cav21c} R. Cavoretto, A. De Rossi, A. Sommariva, M. Vianello, RBFCUB: A numerical package for near-optimal meshless cubature on general polygons, Appl. Math. Lett. 125 (2022), 107704.







\bibitem{dri07} T.A. Driscoll, A.R.H. Heryudono, Adaptive residual subsampling methods for radial basis function interpolation and collocation problems, Comput. Math. Appl. 53 (2007) 927--939.


\bibitem{esm12} M. Esmaeilbeigi, M.M. Hosseini, Dynamic node adaptive strategy for nearly singular problemson large domains, Eng. Anal. Bound. Elem. 36 (2012) 1311--1321.



\bibitem{fas07} G.E. Fasshauer, Meshfree Approximation Methods with \textsc{Matlab}, Interdisciplinary Mathematical Sciences, vol. 6, World Scientific Publishing Co., Singapore, 2007.

\bibitem{fas11} G.E. Fasshauer, Positive definite kernels: Past, present and future, Dolomites Res. Notes Approx. 4 (2011) 21--63.

\bibitem{fas15} G.E. Fasshauer, M.J. McCourt, Kernel-based Approximation Methods using \textsc{Matlab}, Interdisciplinary Mathematical Sciences, Vol. 19, World Scientific Publishing Co., Singapore, 2015.







\bibitem{gao20} K. Gao, G. Mei, S. Cuomo, F. Piccialli, N. Xu, ARBF: adaptive radial basis function interpolation algorithm for irregularly scattered point sets, Soft Computing 24 (2020) 17693--17704.

\bibitem{gao20b} K. Gao, G. Mei, S. Cuomo, F. Piccialli, N. Xu, Adaptive RBF interpolation for estimating missing values in geographical data, in: Y. Sergeyev, D. Kvasov (eds.), Numerical Computations: Theory and Algorithms -- NUMTA 2019, LNCS 11973, pp. 122--130. 



\bibitem{gol15} A. Golbabai, E. Mohebianfar, H. Rabiei, On the new variable shape parameter strategies for radial basis functions, Comput. Appl. Math. 34 (2015) 691--704.





















\bibitem{qia21} B. Qiao, Z. Pan, W. Huang, C. Cao, An adaptive finite-difference method for accurate simulation of first-arrival traveltimes in heterogeneous media, Appl. Math. Comput. 394 (2021) 125792. 








\bibitem{sch95} R. Schaback, Error estimates and condition numbers for radial basis function interpolation, Adv. Comput. Math. 3 (1995) 251--264.



\bibitem{sch11} M. Scheuerer, An alternative procedure for selecting a good value for the parameter $c$ in RBF-interpolation, Adv. Comput. Math. 34 (2011) 105--126. 

\bibitem{sch13} M. Scheuerer, R. Schaback, M. Schlather, Interpolation of spatial data -- A stochastic or a deterministic problem? European J. Appl. Math. 24 (2013) 601--629. 










\bibitem{wen05}  H. Wendland, Scattered Data Approximation, Cambridge Monogr. Appl. Comput. Math., vol. 17, Cambridge Univ. Press, Cambridge, 2005.



\bibitem{zha17} Q. Zhang, Y. Zhao, J. Levesley, Adaptive radial basis function interpolation using an error indicator, Numer. Algorithms 76 (2017) 441--471.

\end{thebibliography}
\end{document}